\documentclass{article}
\usepackage{amsmath}
\usepackage{amssymb}
\begin{document}

\title{ ABOUT LEVI-MALCEV THEOREM FOR HOMOGENEOUS BOL ALGEBRAS}
\author{Thomas B. Bouetou\\ \'Ecole Nationale Sup\'erieure Polytechnique,\\
B.P.~8390 Yaound\'e, Cameroun\\
e-mail:tbouetou@polytech.uninet.cm}

\maketitle \begin{quote}{\bf Abstract}
\em The fundamental ideas of the applicability of Levi-Malcev Theorem for
Bol algebras, which plays a basic role in structural theory are outlined
\end{quote}

  \begin{center}
    {\bf 1. Introduction}
   \end{center}
It is well known that Levi--Malcev theorem is valid for Lie algebras
and Malcev algebras. But it is not valid for binary-Lie algebras (in
classical sense) as well as for Bol algebras. Since the Levi--Malcev
theorem plays the basic role in structural theories, it is useful to
get  some generalized version of such a theorem for Bol algebras. In
this report the problem for so-called homogeneous Bol algebras
with some examples will be discussed.

    Let us note that the problem is important in applications to
Differential geometry and Smooth Quasigroups and Loops theory
(because that Bol algebra is the proper infinitesimal object for a
smooth Bol loop, and any symmetric space, for example, is a smooth Bol
loop).

   It is also noted that, Bol loops and Bol algebras are significant in
Mathematical Physics (chiral anomalies etc.).

Remark:  The smooth Bol loops-Bol algebras theory is due to
 L.V. Sabinin, P.O. Mikheev. See, for example [2].

\begin{center}
{\bf 2. General definitions}
\end{center}

Let consider $({\bf B},\nabla)$ be an affinely connected space with
${\bf R=0}$, and let ${\bf X(x)}$ be an infinitesimal affine
transformation in $({\bf B},\nabla)$ that means for every vector field
${\bf Y} \in B$ we have
 $$ L_X \circ\nabla - \nabla\circ L_X=\nabla_{[X,Y]}$$
where $L$ is a Lie differential in the direction of the vector field
${\bf X}$. Let us introduce
 in ${\bf B}$ a basis of vectori field such that for every
$X_1,\cdot,X_n$ then
 $$X(x)=\alpha^j(x)(X_j)\mid_x $$
and we can have
$$[X_i , X_k]=-T^{j}_{ik}X_j.$$
{\bf Definition 1.} Any vector space $V$ of ${\bf B}$ with the operations
$$\xi,\eta \to \xi\cdot\eta\in V,$$
$$\xi,\eta,\zeta\to (\xi,\eta,\zeta)\in V \quad (\xi,\eta,\zeta\in
V)$$ and identities $$ \xi\cdot\xi=0,\quad
(\xi,\eta,\zeta)+(\eta,\zeta,\xi)+(\zeta,\xi,\eta)=0,$$

$$(\xi,\eta,\zeta)
\chi-(\xi,\eta,\chi)\zeta+(\zeta,\chi,\xi\cdot\eta)-
(\xi,\eta,\zeta\cdot\chi)+\xi\eta\cdot\zeta\chi=0,$$
$$(\xi,\eta,(\zeta,\chi,\omega))=((\xi,\eta,\zeta),\chi,\omega)+(\zeta
,(\xi,\eta,\chi),\omega)+(\zeta,\chi,(\xi,\eta,\omega))$$
is called a Bol algebra.

Let ${\bf I}$ be an ideal of a Bol algebra $V$ and denote
${\bf I}^{(1)}=(V,{\bf I},{\bf I})$ and ${\bf I}^{(k)}=(V,{\bf
I}^{(k-1)},{\bf I}^{(k-1)})$

{\bf Definition 2.} The ideal ${\bf I}$ of a Bol algebra $V$ is called
weakly solvable if there exist k, such that ${\bf I}^{(k)}=0$.


{\bf Definition 3.} The weak radical $RV$ of a Bol algebra is called a
maximal weakly solvable ideal. The Bol algebra $V$ is called semisimple
if $RV=0$. 

{\bf Definition 4.} Let $\cal{ G}$ be a vector subspace of $V$, we say
that $\cal{ G}$ is a Bol subalgebra if, $a,b,c \in {\cal G} \implies
a\cdot b, (a,b,c) \in {\cal G}.$

{\bf Definition 5.} [{\bf Mikheev 1.}] A local analytic loops $(B, \cdot, e)$
verify the $G$-property if there exist such a neighbourhood $\cal{U}$ of the
point $e \in B$ such that for any $x \in \cal{U}$ the loop  $(B, \cdot, e)$ and
 $(B, \frac{1}{x}, e)$ are isomorphic.

{\bf Theorem 1.} [{\bf Mikheev 1.}] A local analytic Bol loop  $(B, \cdot, e)$
posses a  $G$-property if and only, if the corresponding affine connected space
$(B, \Delta)$ is local homogeneous.

{\bf Definition 6.} An arbitrary Bol algebra is said to be homogeneous
if its an infinitesimal object for a local analytic Bol loop that posses a  $G$-property.

\; \; \; \; Indeed we will have the following bilinears and trilinears
operations:
$$
\xi, \eta \longrightarrow \xi \cdot \eta \in T_{e}(B)
$$
$$
(\xi \cdot \eta)^{i}=T^{i}_{jk}(e)\xi^{j}\eta^{k}
$$
$$
(\xi,\eta,\tau)^i=(\Delta_{e}T^{i}_{jk}+T^{s}_{jk}T^{i}_{sl})\xi^{j}\eta^{k}\tau^{l}
$$
$$
[X_i,X_k]=-T^{j}_{ik}X_{j}
$$

which will verify the identities of the definition of Bol algebras.

 Let's $ V $ be a Bol algebra over the field of characteristic zero and let us   assume
 that
in $ V $ there is a radical $ RV $ and a semisimple subalgebra $ SMV $
 such that
$$V=SMV+RV .$$

This implies that $$ RV \cap SMV=\{0\}. $$

Hence $V=RV \oplus SMV $ and $SMV \cong V/RV $
and the Bol subalgebra $ SMV $ is isomorphic to the  factor algebra
of the algebra $ V$ by  its radical.

Connversely, if $V$ contains the Bol subalgebra $SMV$ isomorphic to
$V/RV$ then the algebra $SMV$ is semisimple. Indeed, $SMV \cap RV=0$
and  $$dim V=dim RV+dim (V/RV)= dim RV+dim SMV.$$

Hence, $V=SMV+RV.$
This lead us to state the following theorem:

{\bf Theorem 2.} For every homogenous Bol algebra $V$ over the field of
caracteristic zero there exist semi simple sub algebra $SMV$ and a weak radical
$RV$ such that $V=SMV+RV.$

Proof:

\; \; \; Let assume that $[RV]^2=0$ which means that if the radical is
 abelian, then $RV$ is a submodule of the module $SMV$ relatively to $V$
 by the application $exp({\bf ad}x)$, where ${\bf ad}x$ is an inner
differential of Bol enveloping Lie algebra. Since $[RV]^2=0$ then,
$RV$ belongs to the kernel of the representation of the Bol algebra
$V$ of the defined module $RV.$ Hence we obtain the induced
representation of the algebra $\overline{V}=V/RV.$
 We have for the module the following operation:
$\overline{a\cdot b}=\prod ([a,b]_V))$, $a\in V$, $b\in RV,$ and the
commutator is taken from the minimal enveloping Lie algebra . If
$[RV]^2\neq 0$ then $\overline{V}=(V)/[RV]^2),$ and
$dim\overline{V}\leq dim V.$
 By using the induction according to the dimension of the algebra, we
 can say that the statement is true for the algebra $\overline{V}$.
 Further, $\overline{RV}=
(RV/[RV]^2)$ is a radical for the algebra $\overline{V}$ and
$\overline{V}/\overline{RV}\cong V/RV.$ Therefore, the algebra
$\overline{V}$ contain the subalgebra $\overline{SMV}\cong V/RV$ as a
subalgebra. The algebra $\overline{V}$ and $\overline{SMV}$ has the
form $K/RV$ where $K$ is the subalgebra of Bol algebra
 $V,$ containing $[RV]^2.$ But, $[RV]^2 $ is the radical of the Bol
 algebra
$K$ and $(K)/[RV]^2)\cong (V/RV)$ such that
$dim K \leq dim V,$ using the induction we obtain that $K$ contain the
subalgebra
 $SMV\cong V/RV,$ hence the result.

Further an example is given with a counter example[7].

For the type I of the classification given
in [7], the Bol algebra  have trilinear operation equal to zero and:
$$e_1 \cdot e_3 = e_1 + e_2 $$
$$ e_2 \cdot e_3 = e_2 $$
$$V\cdot V = <e_1+e_2, e_2>$$
$$SMV = \frac{V}{<e_1+e_2, e_2>} = <e_3>.$$

We  see in the case
$$ V = <e_1+e_2, e_2> \oplus <e_3>$$

Counter example, let consider the Bol algebras of Type IV in [7],

$$e_1 \cdot e_3 = xe_1 + pe_2 +e_3, \forall p, x \geq 0 $$
$$ (e_1,e_2, e_3) = e_1 $$
$$ (e_1,e_3, e_3) = e_1. $$
This Bol algebra is not homogeneous.

Hence
$$V\cdot V = <e_1,e_2, e_3>=RV$$
$$SMV = \frac{V}{<e_1,e_2, e_3>} = <e_2,e_3>.$$

We  see that:
$$ RV \cap SMV = <e_1,e_2, e_3> \cap <e_2,e_3>\neq 0.$$

Hence the Levi-Malcev theorem can not be applied.
\begin{center}
{\bf References}
\end{center}

1-P.O. Mikheev {\it On G-properties for analytic Bol loops},Webs and
Quasigroups, Kalinin, 1986. p. 54-59.\\

2-L.V. Sabinin, P.O. Mikheev, {\it Quasigroups and differential
geometry}, in: Quasigroups and loops theory and applications. Berlin:
Helderman Verlas, 1990.-pp. 357-430.\\

3-K. Yamaguti {\it On the Lie triple system and its generalization}
j.sci. Hiroshima Univ. ser. A. 21. N.3 1958. p. 155-160.

4-W.G. Lister {\it A structure theory of Lie triple system} Trans. AMS.
-1952.-T. 72.-p. 217-245.

5-B. Harris {\it Cohomology of Lie triple systems and Lie algebras
with involution} Trans. AMS.-1951.-T. 71.-p. 148-162.

6-M. Kikkawa {\it Remarks on solvability of Lie triple algebras}
Mem. Fac. sci. shim. Univ. N.13 -1979, p.17-22.

7-T.B. Bouetou  {\it On the geometry of Bol algebras} Ph.D Thesis,
Moscow University 1995.   \\
8- V.T. Filipov {\it Homogeneous Bol Algebras} Sibirian Math.
journal Vol.35, N.4 1994, pp. 818-825. \\
9- A.S. Kleshchev {\it Branching Rules for Modular Representation
of Symmetric Groups I } Journal of algebra 178, pp 493-511
(1995).\\

\end{document}